\documentclass[12pt,letterpaper]{amsart}
\usepackage[left=1in,right=1in,top=1.2in,bottom=1.2in]{geometry}
\usepackage{times}
\usepackage{amsmath,amssymb}
\usepackage{graphicx}
\usepackage{booktabs}
\usepackage{xcolor}
\usepackage{placeins}
\usepackage{needspace}
\usepackage{float}
\usepackage{subcaption}
\usepackage[hidelinks]{hyperref}
\makeatletter
\renewcommand\subsection{\@startsection{subsection}{2}{\z@}%
  {0.7\baselineskip}{0.3\baselineskip}{\normalfont\itshape}}
\makeatother

\newcommand{\lab}{f_{n,m}}
\newcommand{\ind}[1]{\mathbf{1}_{\{#1\}}}
\newtheorem{theorem}{Theorem}[section]
\newtheorem{corollary}{Corollary}[section]
\newcommand{\arxiv}[1]{\href{https://arxiv.org/abs/#1}{\texttt{arXiv:#1}}}
\newcommand{\doi}[1]{\href{https://doi.org/#1}{\texttt{doi:#1}}}
\hypersetup{pdftitle={On graceful labeling of Jahangir graphs},
  pdfauthor={Maurizio Imbesi and Carlo Pedala}}
\graphicspath{{figures/}}
\newlength{\constructionpanelgap}
\newlength{\constructionrowgap}
\title[Graceful labeling of Jahangir graphs]{On graceful labeling of Jahangir graphs}
\author[M. Imbesi]{Maurizio Imbesi}
\address[Maurizio Imbesi]{Department of Mathematical and Computer Sciences,
  Physical Sciences and Earth Sciences, University of Messina,
  Viale F. Stagno d'Alcontres 31, I-98166 Messina, Italy}
\email{maurizio.imbesi@unime.it}
\author[C. Pedal\`a]{Carlo Pedal\`a}
\address[Carlo Pedal\`a]{Department of Mathematical and Computer Sciences,
  Physical Sciences and Earth Sciences, University of Messina,
  Viale F. Stagno d'Alcontres 31, I-98166 Messina, Italy}
\email{carlo.pedala@studenti.unime.it}
\date{}

\begin{document}

\begin{abstract}
For a graph with $q$ edges, a graceful labeling assigns distinct
integers from $\{0,\ldots,q\}$ to its vertices so that the absolute
differences across the edges are exactly $1,\ldots,q$.
A Jahangir graph $J_{n,m}$, with $n\geq2$ and $m\geq3$, consists of an
outer cycle on $nm$ vertices and a hub joined to $m$ equally spaced
vertices of that cycle, with $q=(n+1)m$ edges.
This family of graphs is known to be graceful.
In this paper we strengthen the gracefulness property of such graphs by
assigning labels to their vertices through a parametric formula.
These labeled graphs have a monotonicity property: deleting an edge
from the outer cycle leaves a path whose edge differences, in the
negative orientation, are strictly increasing for even $n$
and strictly decreasing for odd $n$.
The starting observation is that large differences admit few endpoint
pairs and require adjacencies between high and low labels.
This suggests a construction based on alternating two separated label
sets along the outer cycle.
The formula has two parity branches, each with two cases modulo $4$.
Applications are given for every case.
\end{abstract}

\maketitle

\section{Introduction}
\label{sec:introduction}

Graph labeling encompasses several types of numerical assignments on graphs.
A labeling assigns numbers to vertices, edges, or both, subject to
prescribed conditions \cite{LopezMuntaner2017}.
A graceful labeling of a graph with $q$ edges assigns distinct integers
from $\{0,\ldots,q\}$ to its vertices and realizes the absolute
differences $1,\ldots,q$ on adjacent vertices.
Rosa introduced these assignments as $\beta$-valuations in 1967
\cite{Rosa1967}; Gardner attributed the term graceful to Golomb
\cite{Gardner1972}. Gallian surveys graph labelings and known graceful
families \cite{Gallian2025}.

The original motivation was cyclic decomposition of complete graphs.
For a graph with $q$ edges, translating graceful vertex labels modulo
$2q+1$ gives isomorphic copies whose edge sets partition $K_{2q+1}$
\cite[Theorem~7]{Rosa1967}; see also
\cite[Theorem~3.7]{ElZanatiVandenEynden2009}.

We use the standard notation $P_t$, $C_t$, $K_t$ for paths, cycles and
complete graphs on $t$ vertices, and $K_{r,s}$ for complete bipartite
graphs with parts of orders $r$ and $s$ \cite{BondyMurty2008}.
Paths admit graceful labelings by alternating small and large values
\cite[Theorem~2]{Rosa1967}.
Caterpillars are trees whose non-leaf vertices form a path.
We write $P_k(a_1,\ldots,a_k)$ for $P_k$ with $a_i$ leaves added at
vertex $i$. Caterpillars are graceful \cite[Theorem~2]{Rosa1967};
further graceful trees and construction operations are surveyed in
\cite{EdwardsHoward2006}.
A cycle $C_t$ is graceful exactly when $t\equiv0$ or $3\pmod4$
\cite[Theorem~1(b)]{Rosa1967}.
Identifying an endpoint of $P_s$ with a vertex of $C_t$ gives the kite
$K_s(C_t)$, which is graceful for $s\geq2$ and $t\geq3$
\cite{HaviarKotulova2024}.
An explicit labeling of $K_3(C_5)$ is given in \cite[Theorem~3.6]{HaviarIvan2026}.

The complete graph $K_t$ is graceful exactly when $t\leq4$;
Golomb proves the obstruction for $t>4$ \cite[Problem~4]{Golomb2008}.
Every complete bipartite graph $K_{r,s}$ is graceful
\cite[Theorem~6]{Rosa1967}.
One may label its two parts by $0,\ldots,r-1$ and $r,2r,\ldots,sr$:
their differences cover $1,\ldots,rs$ exactly once.

Families with a hub lead to Jahangir graphs. The fan $F_t=P_t+K_1$
joins a hub to every path vertex; the wheel $W_t=C_t+K_1$ uses an outer cycle.
Fans are graceful for $t\geq2$
\cite[\S~9, Remarks~(ii)]{GrahamSloane1980}.
Identifying a leaf of a star on $k$ vertices with each path vertex
produces the $k$-enriched fan $kF_t$
\cite[Definition~5.1]{HaviarKurtulik2021}.
These graphs are graceful for $k,t\geq2$
\cite[Theorem~5.3]{HaviarKurtulik2021}.
For wheels, Imbesi and La Barbiera give explicit graceful numberings
\cite[Theorem~3.2]{ImbesiLaBarbiera2024}, with separate
assignments for even and odd outer cycle orders.

\begin{figure}[htbp]
\centering
\includegraphics[width=\linewidth]{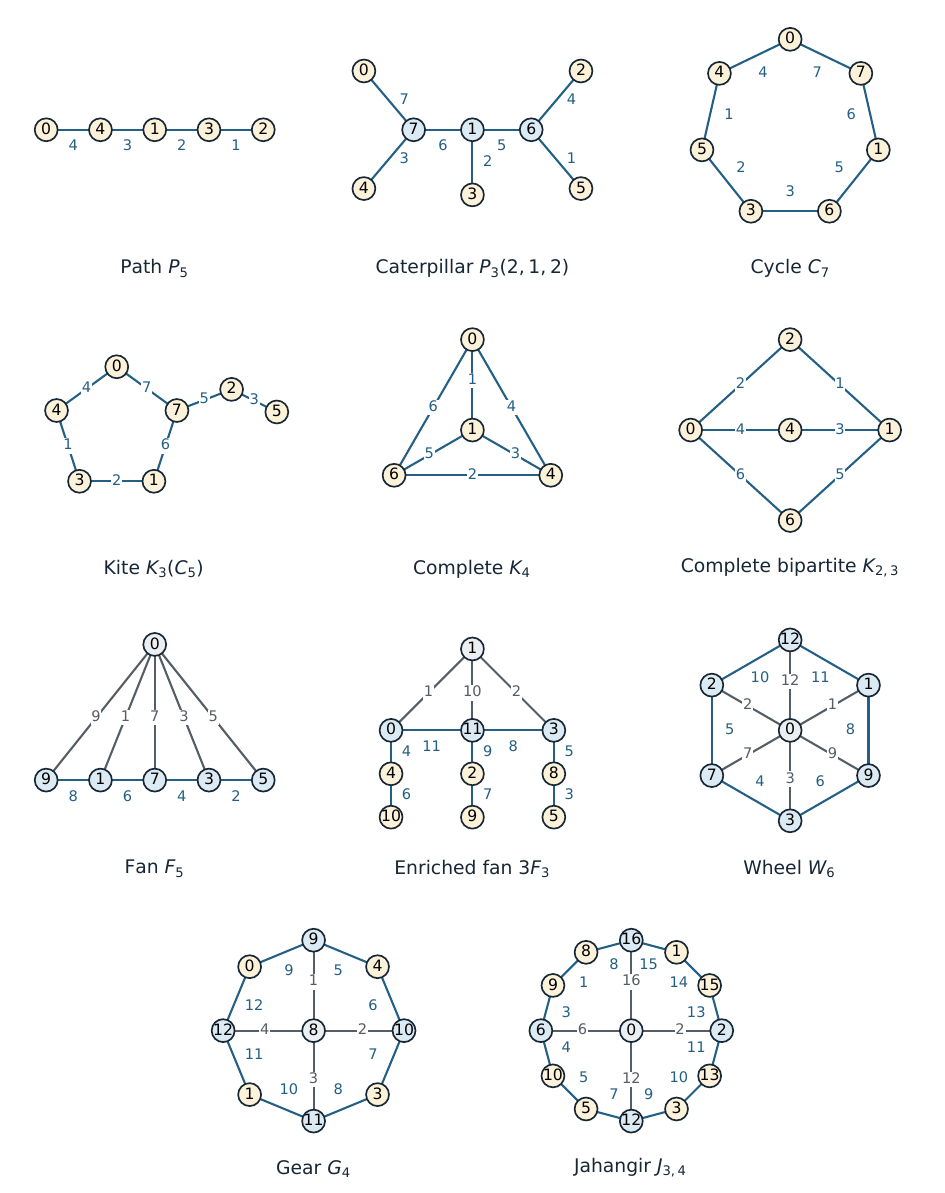}
\caption{Graceful labelings of representative graphs.}
\label{fig:graceful-families}
\end{figure}

The gear $G_m$ subdivides each outer cycle edge of $W_m$ once, retaining
the spokes. Replacing each outer cycle edge by a path of $n$ edges gives
the multigear $(nG)_m$
\cite[Definition~3.7]{ImbesiLaBarbiera2024}.
The Jahangir graph $J_{n,m}$ consists of an outer cycle on $nm$ vertices and
a hub joined to $m$ equally spaced outer cycle vertices
\cite[Definition~3.9]{ImbesiLaBarbiera2024}.
Thus $n$ counts outer cycle edges between consecutive radial vertices and
$m$ counts spokes. We have $W_m=J_{1,m}$
\cite[Remark~3.10]{ImbesiLaBarbiera2024}, while
$G_m=J_{2,m}$ and $(nG)_m=J_{n,m}$
\cite[Remark~3.11]{ImbesiLaBarbiera2024}.
These graphs are graceful \cite[Proposition~3.12]{ImbesiLaBarbiera2024};
explicit gear labelings are given in
\cite[Theorem~3.14]{ImbesiLaBarbiera2024}.
Representative graceful labelings appear in Figure~\ref{fig:graceful-families}.

Rofa gives a direct algebraic labeling
function for rooted symmetric trees \cite[Theorem~1]{Rofa2021}.
For Jahangir graphs, we give an explicit graceful labeling
for every $n\geq2$ and $m\geq3$, with monotone differences along the outer cycle after removing one edge.
The repeated sectors give a natural way to describe vertex positions,
but the outer cycle and spokes must be treated together: radial vertex
labels occur in both sets of edge differences.

The largest differences are especially restrictive. The maximum
requires adjacency between zero and the maximum available label,
while the next difference permits only two endpoint pairs.
This suggests reserving differences for the spokes and closing edge,
then organizing the remaining assignments along the path left by removing
that edge.
Our construction alternates high and low label bands along this path
and orders its differences monotonically.
The restriction on the largest differences is necessary for every
graceful labeling; the two bands and monotone order are choices of
the construction. The hub and radial labels are chosen jointly to
realize distinct spoke differences and keep the hub value unused
on the outer cycle. Deleting one specified outer cycle edge leaves a path whose
differences increase for even $n$ and decrease for odd $n$ in the
negative orientation. The formula has two parity branches, each
with two cases modulo $4$.

In Section~\ref{sec:definitions}, we fix the vertex notation and the
orientation of the outer cycle, then discuss the observations that motivate
the construction. In Section~\ref{sec:construction}, we state the theorem
and prove it by constructing the labeling. The proof treats even and odd
$n$ separately, giving a common construction and verification for the cases
of each parity. The corollary follows from this construction. Finally,
Section~\ref{sec:examples} shows how to use the formula on individual
Jahangir graphs. Each application presents the assignments in steps,
relating the parameters of the formula to the vertex labels and the induced
edge labels.

\FloatBarrier
\section{Definitions and motivation}
\label{sec:definitions}
\begingroup
\setlength{\intextsep}{4pt}

For integers $n\geq2$ and $m\geq3$, the Jahangir graph $J_{n,m}$
consists of an outer cycle of $N=nm$ vertices and a hub $h$ joined to $m$
equally spaced outer cycle vertices \cite[Definition~3.9]{ImbesiLaBarbiera2024}.
It has $N+1$ vertices and
$q=(n+1)m$ edges. Write the outer cycle vertices as $v_{i,t}$ in sector
order, with $0\leq i<m$ and $0\leq t<n$. The radial vertices
$v_{i,0}$ are adjacent to the hub. Sector $i$ is the path
\[
v_{i,0},v_{i,1},\ldots,v_{i,n-1},v_{(i+1)\bmod m,0},
\]
of $n$ outer cycle edges; its terminal vertex is indexed in the next
sector. We take the negative orientation of the outer cycle,
following the sector order above.
Figure~\ref{fig:structure} shows the sectors of $J_{3,3}$.

\begin{figure}[H]
\centering
\includegraphics[width=\linewidth]{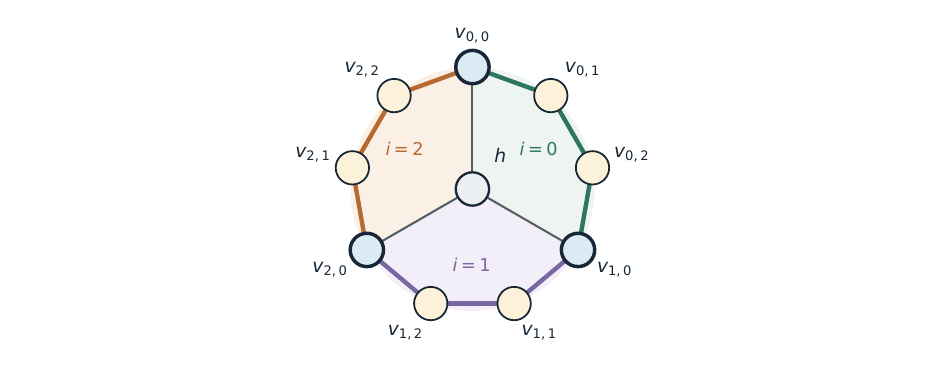}
\caption{The graph $J_{3,3}$ and its sectors.}
\label{fig:structure}
\end{figure}

A graceful labeling $f$ of a graph with $q$ edges assigns distinct
integers from $\{0,\ldots,q\}$ to its vertices so that the absolute
edge differences are exactly $1,\ldots,q$ \cite{Rosa1967}.
The edge labels are induced by the vertex labels: each edge $uv$
receives the value $|f(u)-f(v)|$.
Figure~\ref{fig:graceful-example}
shows such an assignment for $J_{3,3}$.

\begin{figure}[H]
\centering
\includegraphics[width=\linewidth]{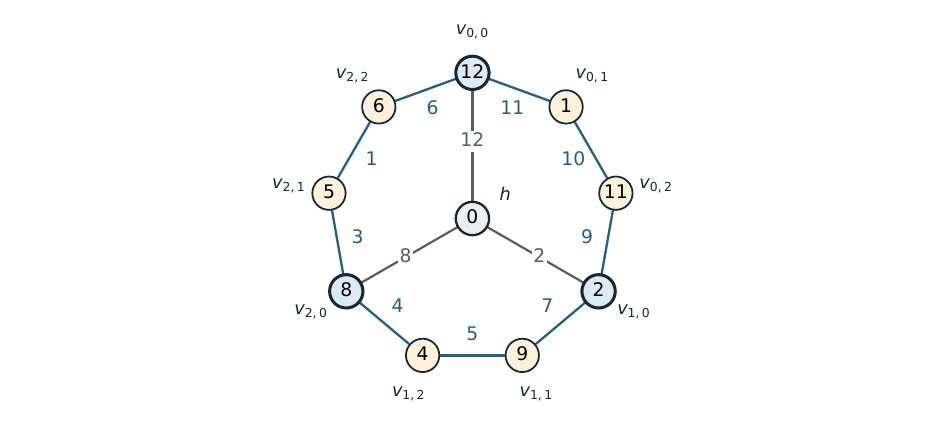}
\caption{A graceful labeling of $J_{3,3}$.}
\label{fig:graceful-example}
\end{figure}

The largest differences admit few endpoint pairs.
Difference $q$ requires the pair $(0,q)$; difference $q-1$ allows only
$(0,q-1)$ or $(1,q)$. Decreasing the difference by one adds one pair.
In general, difference $d$, with $1\leq d\leq q$, requires labels
$(z,z+d)$ with $0\leq z\leq q-d$, giving $q-d+1$ pairs:
\[
\begin{array}{c|l|c}
\text{Difference}&\text{Possible endpoint pairs}&\text{Number}\\ \hline
q&(0,q)&1\\
q-1&(0,q-1),\ (1,q)&2\\
q-2&(0,q-2),\ (1,q-1),\ (2,q)&3\\
d&(z,z+d),\quad 0\leq z\leq q-d&q-d+1.
\end{array}
\]
Figure~\ref{fig:candidate-pairs} illustrates the count for $q=12$.

\begin{figure}[H]
\centering
\includegraphics[width=\linewidth]{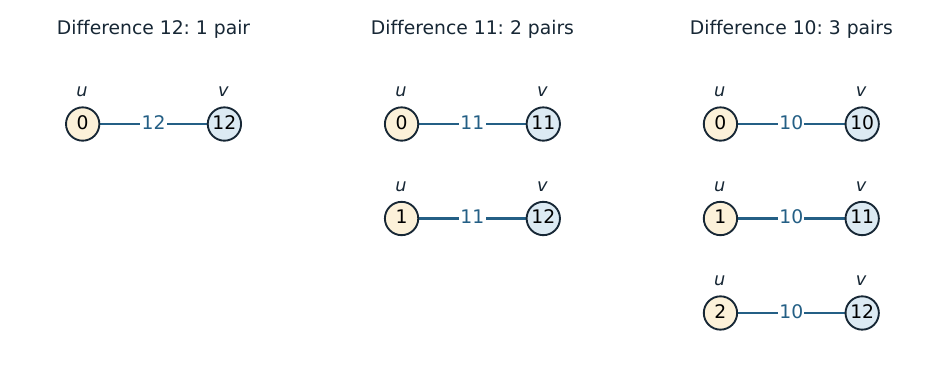}
\caption{Candidate endpoint pairs for $q=12$.}
\label{fig:candidate-pairs}
\end{figure}

This restriction forces high--low adjacencies for the largest
differences. If $d>q/2$, its endpoint labels satisfy
\[
 z\leq q-d<q/2,\qquad z+d\geq d>q/2.
\]
The endpoints lie on opposite sides of $q/2$.
Every graceful labeling must realize all
$\lceil q/2\rceil$ such differences. Since only $m$ edges are spokes,
at least $\lceil q/2\rceil-m$ of these high--low adjacencies lie on
the outer cycle, as illustrated by the pairs in Figure~\ref{fig:candidate-pairs}.

These forced adjacencies suggested two separated label bands, $H$
and $L$, with $\max L<\min H$. We alternate them along the outer cycle,
apart from the closing edge when $n$ and $m$ are both odd: that edge joins
two $H$ labels. All outer cycle vertices receive a label
(Figure~\ref{fig:alternation}).

\begin{figure}[H]
\centering
\includegraphics[width=\linewidth]{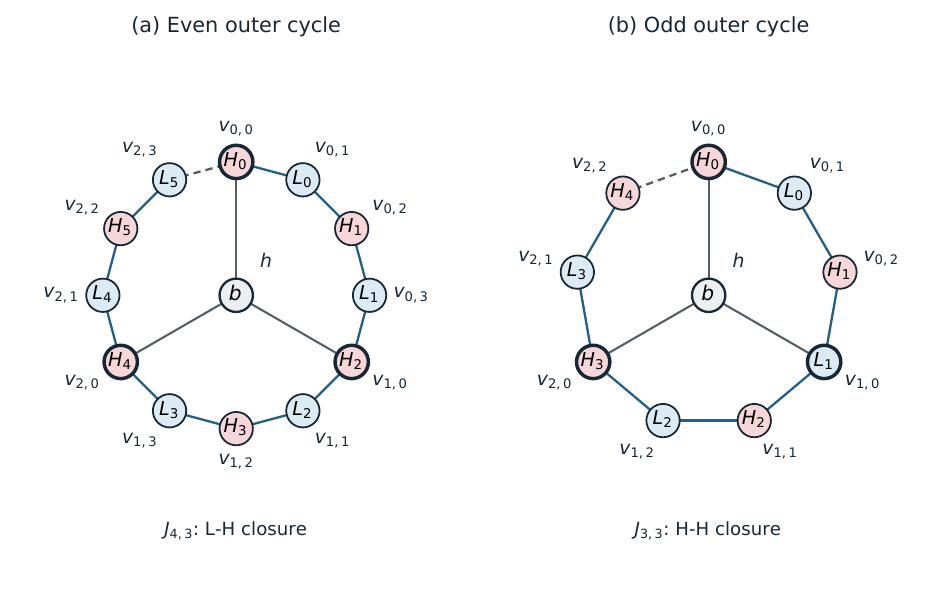}
\caption{High and low label positions for even (a) and odd (b) outer cycles;
$b=\lab(h)$.}
\label{fig:alternation}
\end{figure}

We reserve the spoke differences and a closing difference, then seek
a monotone ordering of the remaining outer cycle differences
(Figure~\ref{fig:monotonicity}). This ordering and the alternation of
bands are choices of the construction; neither is forced by gracefulness.

\begin{figure}[H]
\centering
\includegraphics[width=\linewidth]{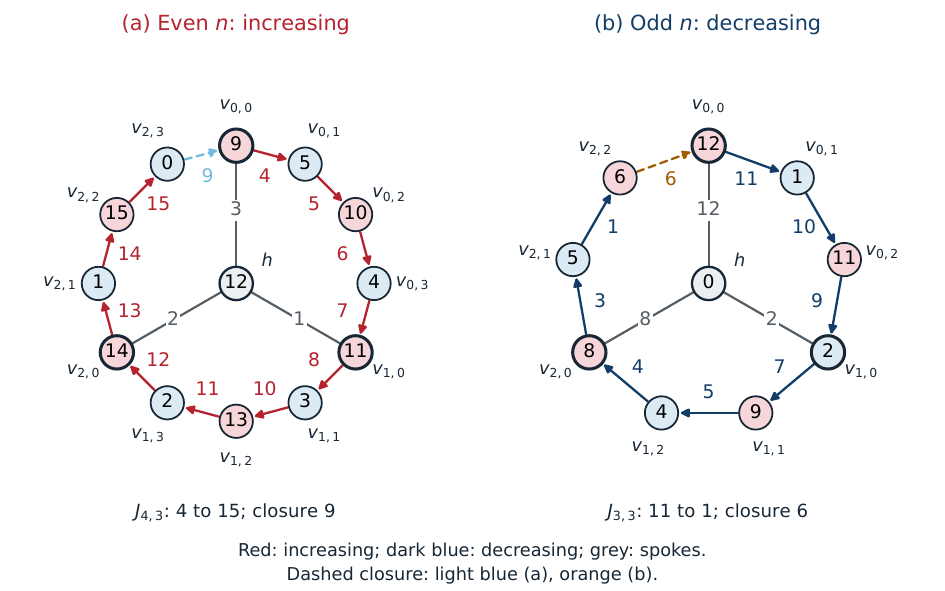}
\caption{Strictly monotone differences along the outer cycle with one edge removed.}
\label{fig:monotonicity}
\end{figure}

The hub and radial assignments must realize the reserved spoke
differences while keeping the hub label unused on the outer cycle.
Spoke differences need not be monotone.
\endgroup

\Needspace{12\baselineskip}
\section{Main result}
\label{sec:construction}

\begin{theorem}\label{thm:main}
If $J_{n,m}$ is a Jahangir graph with $n\geq2$ and $m\geq3$, then
there exists a function $f_{n,m}$ which defines a graceful labeling of
$J_{n,m}$: the vertex labels are distinct integers from $\{0,\ldots,q\}$,
and the induced edge labels take each value in $\{1,\ldots,q\}$ exactly
once, where $q=(n+1)m$.
\end{theorem}

\begin{proof}
We construct $f_{n,m}$ in four cases according to $n\bmod4$.
All indices and intervals below are integer-valued; empty intervals
are omitted. Write $\ind{P}=1$ when $P$ holds and $0$ otherwise.
The closing edge joins the endpoints of the auxiliary path and
has difference $c$.
The formulas below align the path vertices with the sectors and
assign the hub value.

\medskip\noindent\textbf{$n$ even.}\par\smallskip

Let $n=2k$, $P=km$ and $a=(k+1)m$.
An even sector has endpoints in the same label band; radial vertices
are $k$ high positions apart. The parity of $k$ distinguishes
$n\equiv0,2\pmod4$.
Each even case specifies the hub label $\lab(h)=b$, an increasing sequence
$H_0,\ldots,H_{P-1}$, a phase $\eta$ and $c=H_0$.
Reserve the spoke differences
\[
 S=\{|H_{ki+\eta}-b|:0\leq i<m\},
\]
and let $d_1<\cdots<d_{2P-2}$ be the elements of
$\{1,\ldots,q\}\setminus(S\cup\{c,q\})$. Set
\[
 L_j=H_j-d_{2j+1}\quad(0\leq j<P-1),
 \qquad L_{P-1}=0.
\]
The labels are
\begin{equation}\label{eq:even-labels}
 \begin{aligned}
 \lab(v_{i,2s})&=H_{(ki+s+\eta)\bmod P},\\
 \lab(v_{i,2s+1})&=L_{(ki+s+\eta)\bmod P},
 \qquad 0\leq i<m,\quad 0\leq s<k .
 \end{aligned}
\end{equation}
Thus the outer cycle is the list $H_0,L_0,\ldots,H_{P-1},L_{P-1}$
rotated to start at $H_\eta$.
The high sequence fixes alternate outer cycle positions, including all
radial vertices. Their distances from the hub give $S$; the residual
differences then determine the low labels.

Each segment $H_j,L_j,H_{j+1}$ realizes a pair of differences through
the identity
\begin{equation}\label{eq:pairs}
 H_{j+1}-H_j=d_{2j+2}-d_{2j+1}
 \quad\Longrightarrow\quad
 \begin{cases}
 H_j-L_j=d_{2j+1},\\
 H_{j+1}-L_j=d_{2j+2}.
 \end{cases}
\end{equation}
Once their gaps match the high increments, every residual difference
is realized. Moreover,
\[
 L_j-L_{j+1}=d_{2j+3}-d_{2j+2}>0
 \qquad(0\leq j<P-2).
\]
In the cases below, $H_{P-1}=q$, $H_{P-2}=q-1$,
the last residual pair is $(q-2,q-1)$, and $L_0<a\leq H_0$.
Consequently $L_{P-2}=1$, and the low labels decrease strictly to zero.
The two bands are disjoint and lie in $[0,q]$.
The two edges incident with $L_{P-1}=0$ supply $q$ and $c$.

\Needspace{3\baselineskip}
\medskip\noindent\textit{Case $n\equiv0\pmod4$.}\par\smallskip

Here $k$ is even. Let $g=P/2$ and set $\lab(h)=b=a+g$. Put
\[
 c=a,\qquad
 \eta=\ind{m\text{ odd and }k\equiv0\pmod4},
\]
and define
\begin{equation}\label{eq:case-zero}
 H_j=a+j+\ind{j\geq g}\qquad(0\leq j<P).
\end{equation}
For $m=2s$ the reserved spoke set is
\[
S=\{1+kj,k+kj:0\leq j<s\}.
\]
For $m=2s+1$, put $t=k/2-\eta$. Then
\[
S=\{t+kj:0\leq j\leq s\}
 \cup\{k-t+1+kj:0\leq j<s\}.
\]

\medskip\noindent\textit{Case $n\equiv2\pmod4$.}\par\smallskip

Here $k=n/2$ is odd, including $k=1$. Set
\[
 \lab(h)=b=
 \begin{cases}
 a,&m\text{ even},\\
 a-k-1,&m\text{ odd}.
 \end{cases}
\]
Take $\eta=0$ and use the following outer cycle parameters in
\eqref{eq:even-labels}:
\begin{equation}\label{eq:case-two}
 \begin{array}{c|cc}
 &c&H_j\quad(0\leq j<P)\\ \hline
 m\text{ even}&a+1&a+1+j\\[2pt]
 m\text{ odd}&a&a+j+\ind{j\geq(k+1)/2}.
 \end{array}
\end{equation}
The spoke sets are respectively
\[
 S=\{1+ki:0\leq i<m\},
 \qquad
 S=\{k+1\}\cup\{k(i+1)+2:1\leq i<m\}.
\]

\Needspace{7\baselineskip}
\medskip\noindent\textbf{Verification for both even cases.}\par\smallskip
For these spoke sets and $0\leq j<P-1$, the retained intervals give
\[
d_{2j+2}-d_{2j+1}=H_{j+1}-H_j
=1+\begin{cases}
\ind{j=P/2-1},& k\text{ even},\\
\ind{j=(k-1)/2},& k,m\text{ odd},\\
0,&k\text{ odd},\ m\text{ even}.
\end{cases}
\]
Indeed consecutive omissions leave intervals of even length,
apart from the single pair crossing the displayed jump of $H$.
Identity~\eqref{eq:pairs} therefore realizes every difference in the
residual list. The two remaining outer cycle edges have differences $q$ and $c$.

The hub value is omitted from the high sequence when $k$ is even,
and lies between the two bands when $k$ is odd and $m$ is even.
When $k,m$ are odd, the low sequence satisfies
\[
L_{k-1}=b+1,\qquad L_k\leq b-1,
\]
so it skips the hub as well. Thus all vertex labels are distinct
and lie in $[0,q]$. The spoke sets contain $m$ distinct positive
values, disjoint from $\{c,q\}$; together with the residual list
and the two outer cycle edges at $L_{P-1}$, they partition $\{1,\ldots,q\}$.

\medskip\noindent\textbf{$n$ odd.}\par\smallskip

Let $n=2p-1$, $\varepsilon=p\bmod2$, $K=\lceil m/2\rceil$ and
$L=\lfloor m/2\rfloor$.
An odd sector contained in the open path has endpoints in opposite
label bands. The parity of $p$ determines the phase, distinguishing
$n\equiv1,3\pmod4$. Both cases use
\[
\begin{aligned}
S={}&\{q-\varepsilon-2pr:0\leq r<K\}
     \cup\{p+2pr:0\leq r<L\},\\
c={}&pm+\varepsilon\ind{m\text{ even}}.
\end{aligned}
\]
Order $\{1,\ldots,q\}\setminus(S\cup\{c\})$ decreasingly as
$d_1,\ldots,d_{N-1}$ and put
\begin{equation}\label{eq:case-one-path}
x_0=q,\qquad x_j=q+\sum_{r=1}^{j}(-1)^r d_r
\quad(1\leq j<N).
\end{equation}
The assignment on the outer cycle is
\begin{equation}\label{eq:case-one}
\lab(v_{i,t})=x_{(ni+t+\varepsilon)\bmod N}
\quad(0\leq i<m,\ 0\leq t<n).
\end{equation}
The absolute differences along the auxiliary path are exactly
$d_1,\ldots,d_{N-1}$.

\medskip\noindent\textit{Case $n\equiv1\pmod4$.}\par\smallskip
Here $p\geq3$ is odd and $\varepsilon=1$. Assign
\[
\lab(h)=q-1.
\]
Equations~\eqref{eq:case-one-path} and~\eqref{eq:case-one} give
all outer cycle labels.

\medskip\noindent\textit{Case $n\equiv3\pmod4$.}\par\smallskip
Here $p$ is even and $\varepsilon=0$. Assign
\[
\lab(h)=0.
\]
The same outer cycle assignment has the following sector form. Put
\[
r=p/2,\qquad M=\lfloor m/2\rfloor,\qquad
A_i=pi,\qquad B_i=q-pi,\qquad\delta_i=\ind{i\geq M}.
\]
Then
\begin{equation}\label{eq:case-three}
\begin{array}{c|cc}
&\lab(v_{i,2s})&\lab(v_{i,2s+1})\\ \hline
i\text{ even}&B_i-s-\delta_i\ind{s\geq r}&A_i+s+1\\[2pt]
i\text{ odd}&A_i+s+\delta_i\ind{s\geq r}&B_i-s-1.
\end{array}
\end{equation}
Here $0\leq i<m$, with $0\leq s<p$ in even positions and
$0\leq s<p-1$ in odd positions. To identify this expression with
\eqref{eq:case-one-path}, put $t_i=q-2pi$. For $0\leq i<m-1$,
the differences through a sector and the next radial vertex are
\[
t_i-j\quad(1\leq j<p),\qquad
t_i-j-\delta_i\quad(p\leq j<2p).
\]
The last sector has the internal difference set
$\{1,\ldots,2p-1\}\setminus\{p\}$. These are precisely the
successive blocks of the descending residual list. Since the first
label is $q$, the two expressions agree.

\Needspace{7\baselineskip}
\medskip\noindent\textbf{Verification for both odd cases.}\par\smallskip
For both odd cases the hub value is $b=\varepsilon(q-1)$.
The smallest element of the upper progression in $S$ and the
largest element of the lower progression are
\[
U=q-\varepsilon-2p(K-1),\qquad V=p+2p(L-1).
\]
They satisfy $U-V=3p-\varepsilon$ and $V<c<U$; hence $S$ has
$m$ distinct values and does not contain $c$.

For $0\leq i<m$, set $T_i=q-2pi-\varepsilon-1$.
The numbers of reserved values above $T_i$ in the upper progression,
the lower progression and $\{c\}$ are respectively
$\min(i+1,K)$, $\max(0,i-K)$ and $\ind{i\geq K}$.
Their sum is $i+1$, so $T_i=d_{ni+\varepsilon+1}$.
For $0\leq i<m-1$, the integer interval
$[T_i-2p+1,T_i]$ therefore contains exactly one reserved value.
The other $2p-1$ values form the block between successive radial
positions. All reserved values have parity $\varepsilon$, opposite
to $T_i$, so after its first term this block consists of $p-1$
consecutive pairs. Its alternating sum, with positive first sign, is
$T_i-(p-1)=q-(2i+1)p-\varepsilon$. Hence
\[
x_{n(i+1)+\varepsilon}-x_{ni+\varepsilon}
=(-1)^{i+\varepsilon+1}\bigl(q-(2i+1)p-\varepsilon\bigr)
\quad(0\leq i<m-1).
\]
Starting from $x_{\varepsilon}=(1-\varepsilon)q$, induction yields
the radial labels
\begin{equation}\label{eq:case-one-radials}
R_i=\lab(v_{i,0})=
\begin{cases}
q-pi-\varepsilon,&i+\varepsilon\text{ even},\\
pi,&i+\varepsilon\text{ odd}.
\end{cases}
\end{equation}
Their distances from $b$ give exactly $S$. The final retained block
consists of $\{1,\ldots,2p-\varepsilon-1\}\setminus\{p\}$ in
decreasing order. Its alternating sum with positive first sign is
$p$: writing $p=2r$ gives $3r-r=2r$, while $p=2r+1$ gives
$(3r+1)-r=2r+1$. Thus
\[
x_{N-1}=R_{m-1}+(-1)^{m+\varepsilon}p=q-c.
\]
Thus the edge joining the path endpoints has difference $c$.

Finally, the strict decrease of the $d_j$ makes $x_0,x_2,\ldots$
strictly decreasing and $x_1,x_3,\ldots$ strictly increasing.
Their closest terms are consecutive along the path and differ by
$d_{N-1}=1$. The bands are therefore disjoint and lie in $[0,q]$.
When $\varepsilon=0$ all path labels are positive, leaving the hub
value $0$ unused. When $\varepsilon=1$, the first labels are
$q,0,q-2$, leaving the hub value $q-1$ unused. The path, closure
and spokes realize the disjoint sets
$\{d_1,\ldots,d_{N-1}\}$, $\{c\}$ and $S$, respectively.
This proves gracefulness in all four cases.
\end{proof}

\Needspace{8\baselineskip}
\begin{corollary}\label{cor:monotonicity}
For the labeling constructed in the proof of Theorem~\ref{thm:main},
deleting the closing edge from the outer cycle leaves a path whose
edge differences, in the negative orientation, are strictly increasing
if $n$ is even and strictly decreasing if $n$ is odd.
\end{corollary}

\begin{proof}
For even $n$, identity~\eqref{eq:pairs} gives the path differences
$d_1,\ldots,d_{2P-2},q$ in this order, with
$d_1<\cdots<d_{2P-2}<q$.
For odd $n$, \eqref{eq:case-one-path} gives the strictly decreasing
list $d_1,\ldots,d_{N-1}$.
In both cases the omitted edge has difference $c$.
\end{proof}

\section{Applications}
\label{sec:examples}

We apply the formulas with $m=4$ and $n=4,5,6,7$, one graph per
residue class (Table~\ref{tab:examples}). Each assignment step has
a corresponding panel.

\begin{table}[htbp]
\centering
\caption{Parameters of the four applications; $S$ contains the spoke differences.}
\label{tab:examples}
\begin{tabular}{ccrrl}
\toprule
$n\bmod4$ & Graph & $q$ & $\lab(h)$ & $S$\\
\midrule
0 & $J_{4,4}$ & 20 & 16 & $\{1,2,3,4\}$\\
1 & $J_{5,4}$ & 24 & 23 & $\{3,9,17,23\}$\\
2 & $J_{6,4}$ & 28 & 16 & $\{1,4,7,10\}$\\
3 & $J_{7,4}$ & 32 & 0 & $\{4,12,24,32\}$\\
\bottomrule
\end{tabular}
\end{table}
\FloatBarrier
\clearpage
\subsection{\texorpdfstring{$J_{4,4}$}{J(4,4)}}

\begin{figure}[p]
\centering
\makebox[\linewidth][c]{%
\begin{subfigure}[t]{0.332\linewidth}
\centering
\includegraphics[width=\linewidth]{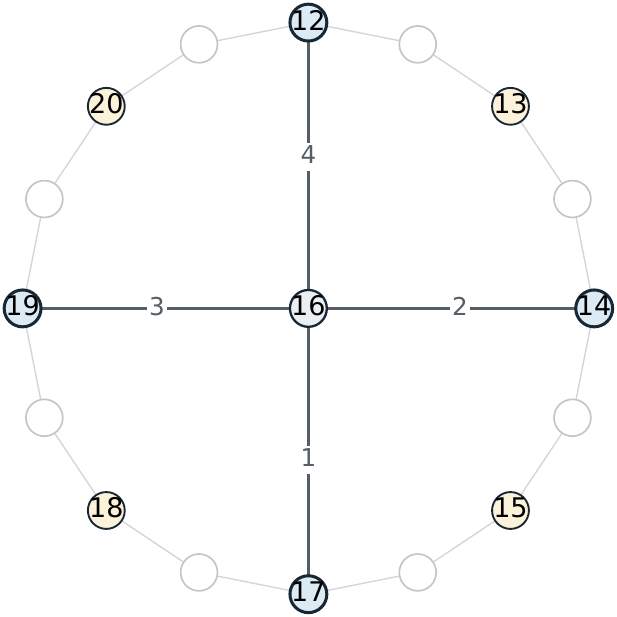}
\caption{Hub and high labels.}
\label{fig:j44-step1}
\end{subfigure}%
\hspace{\constructionpanelgap}%
\begin{subfigure}[t]{0.332\linewidth}
\centering
\includegraphics[width=\linewidth]{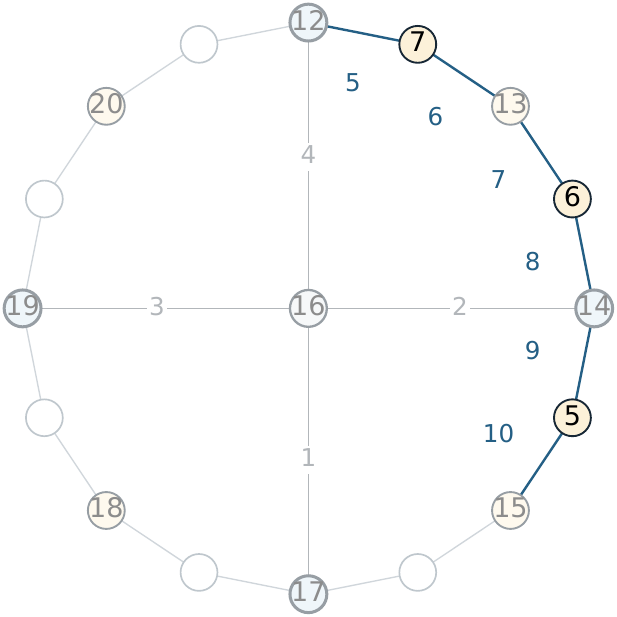}
\caption{First three low labels.}
\label{fig:j44-step2}
\end{subfigure}%
}\par\vspace{\constructionrowgap}
\makebox[\linewidth][c]{%
\begin{subfigure}[t]{0.332\linewidth}
\centering
\includegraphics[width=\linewidth]{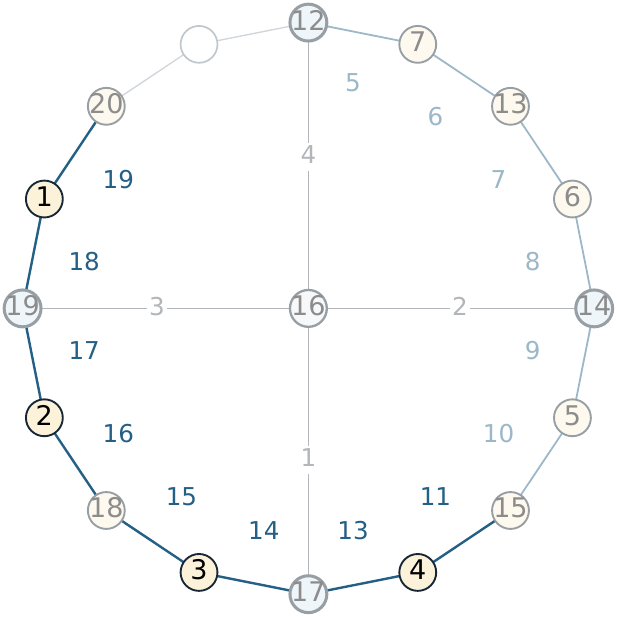}
\caption{Remaining low labels.}
\label{fig:j44-step3}
\end{subfigure}%
\hspace{\constructionpanelgap}%
\begin{subfigure}[t]{0.332\linewidth}
\centering
\includegraphics[width=\linewidth]{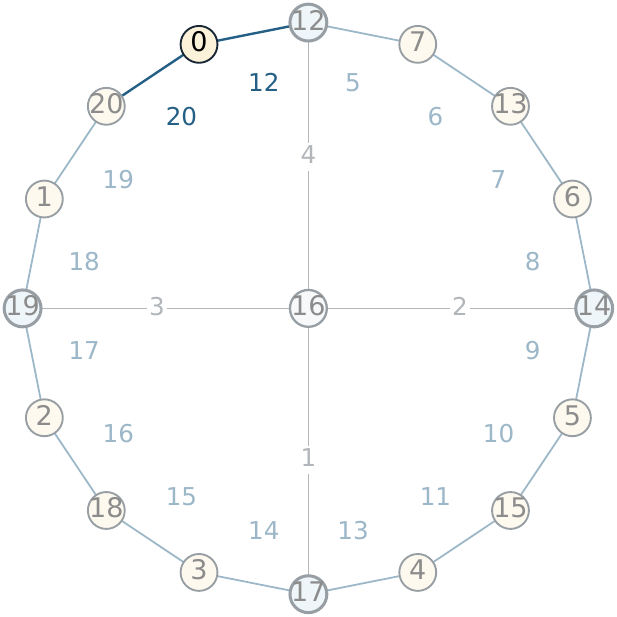}
\caption{Zero and outer cycle closure.}
\label{fig:j44-step4}
\end{subfigure}%
}\par\smallskip
\caption{Construction of $J_{4,4}$.}
\label{fig:j44}
\end{figure}

\noindent\textit{Step 1. Hub and high labels (Figure~\ref{fig:j44-step1}).} 
Here $k=2$, $P=8$, $a=12$, $g=4$ and $q=20$. Equation~\eqref{eq:case-zero} gives
\[
 b=16,\qquad \eta=0,\qquad H_j=12+j+\ind{j\geq4},
 \qquad H=(12,13,14,15,17,18,19,20).
\]
Place these labels at alternate outer cycle positions. Assign the omitted
value $16$ to the hub. The radial labels $H_0=12$, $H_2=14$, $H_4=17$
and $H_6=19$ give spoke differences $4$, $2$, $1$ and $3$.
Reserve $S=\{1,2,3,4\}$ and the differences $c=12$, $q=20$
through the final zero. The residual list is
\[
 D=(5,6,7,8,9,10,11,13,14,15,16,17,18,19).
\]
\par\medskip

\noindent\textit{Step 2. First three low labels (Figure~\ref{fig:j44-step2}).} 
Apply $L_j=H_j-d_{2j+1}$ to $(5,6)$, $(7,8)$ and $(9,10)$:
\[
 L_0=12-5=7,\qquad L_1=13-7=6,\qquad L_2=14-9=5.
\]
The segments $(12,7,13)$, $(13,6,14)$ and $(14,5,15)$ realize
these three pairs of differences.
\par\medskip

\noindent\textit{Step 3. Remaining positive low labels (Figure~\ref{fig:j44-step3}).} 
The pair $(11,13)$ matches the high increment $17-15=2$. It gives
$L_3=15-11=4$. The remaining pairs give
\[
 L_4=17-14=3,\qquad L_5=18-16=2,\qquad L_6=19-18=1.
\]
The segments through $L_3,L_4,L_5,L_6$ realize $(11,13)$,
$(14,15)$, $(16,17)$ and $(18,19)$, respectively. The first pair
allows the high sequence to skip the hub label $16$.
\par\medskip

\noindent\textit{Step 4. Zero and outer cycle closure (Figure~\ref{fig:j44-step4}).} 
Set $L_7=0$. The last two edges have differences
\[
 |20-0|=20,\qquad |0-12|=12.
\]
The completed outer cycle, grouped by sectors, is
\[
 (12,7,13,6)\mid(14,5,15,4)\mid(17,3,18,2)\mid(19,1,20,0).
\]
Its differences are $\{5,\ldots,20\}$; the spokes supply $\{1,2,3,4\}$.
\par\medskip

\clearpage

\subsection{\texorpdfstring{$J_{5,4}$}{J(5,4)}}

\begin{figure}[p]
\centering
\makebox[\linewidth][c]{%
\begin{subfigure}[t]{0.332\linewidth}
\centering
\includegraphics[width=\linewidth]{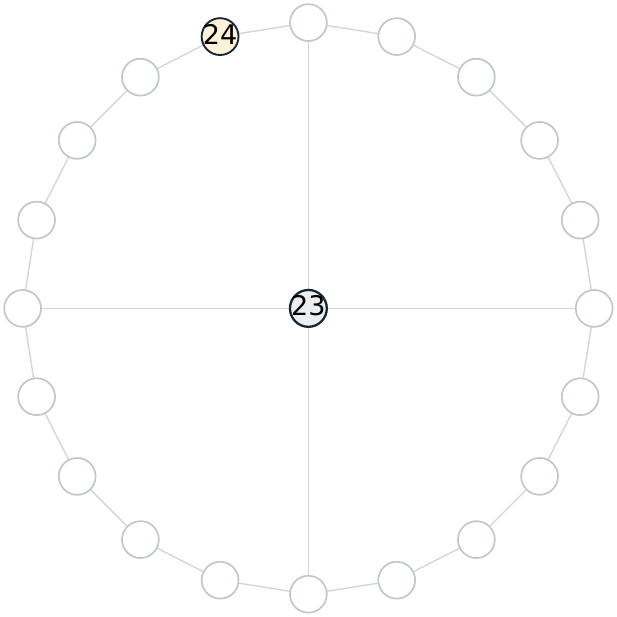}
\caption{Hub and starting label.}
\label{fig:j54-step1}
\end{subfigure}%
\hspace{\constructionpanelgap}%
\begin{subfigure}[t]{0.332\linewidth}
\centering
\includegraphics[width=\linewidth]{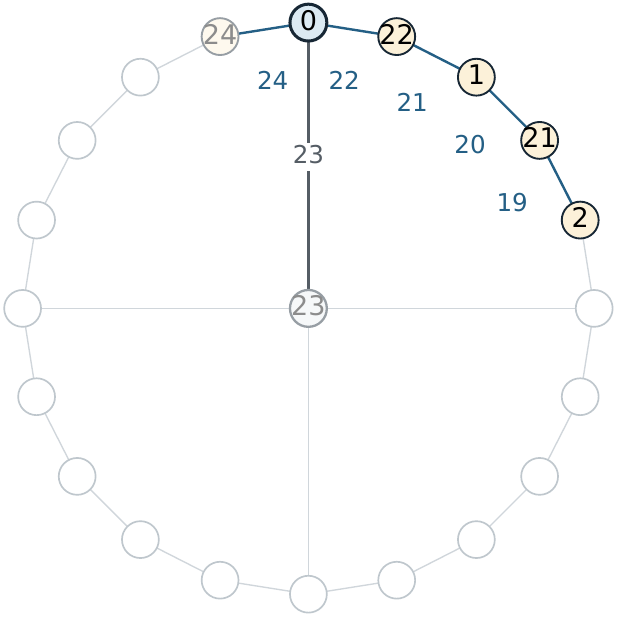}
\caption{First sector.}
\label{fig:j54-step2}
\end{subfigure}%
\hspace{\constructionpanelgap}%
\begin{subfigure}[t]{0.332\linewidth}
\centering
\includegraphics[width=\linewidth]{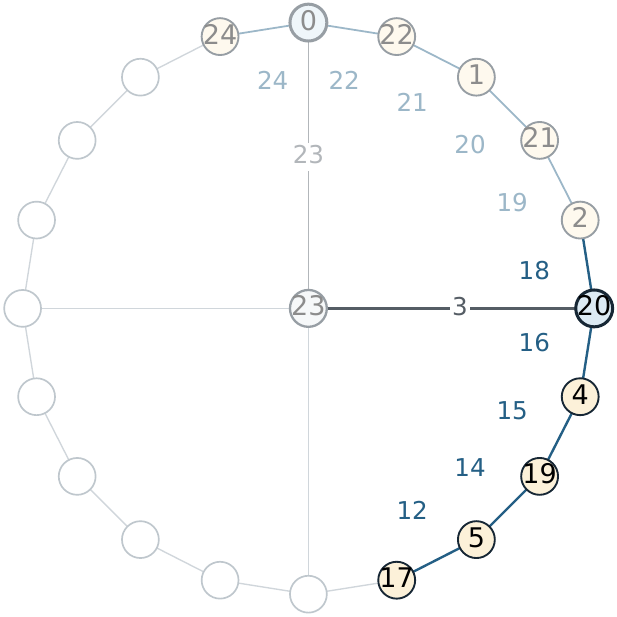}
\caption{Second sector.}
\label{fig:j54-step3}
\end{subfigure}%
}\par\vspace{\constructionrowgap}
\makebox[\linewidth][c]{%
\begin{subfigure}[t]{0.332\linewidth}
\centering
\includegraphics[width=\linewidth]{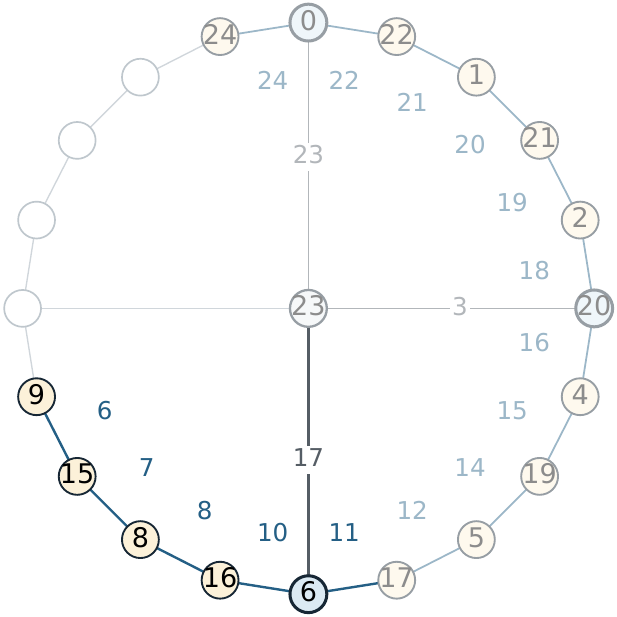}
\caption{Third sector.}
\label{fig:j54-step4}
\end{subfigure}%
\hspace{\constructionpanelgap}%
\begin{subfigure}[t]{0.332\linewidth}
\centering
\includegraphics[width=\linewidth]{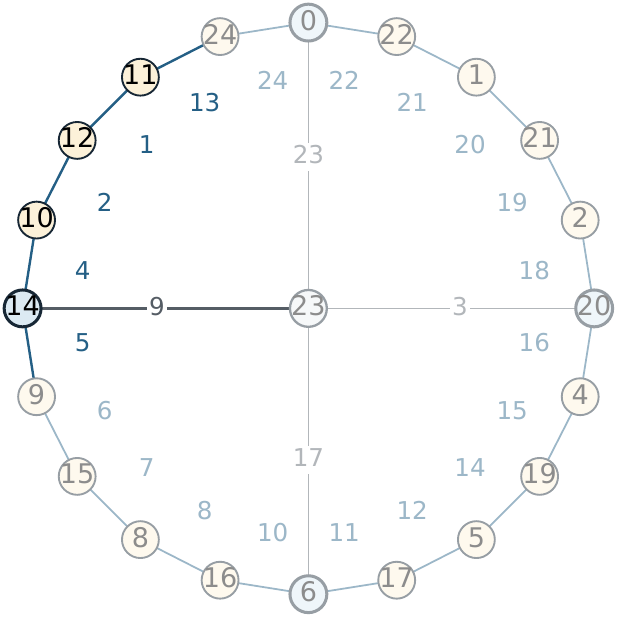}
\caption{Last sector and closure.}
\label{fig:j54-step5}
\end{subfigure}%
}\par\smallskip
\caption{Construction of $J_{5,4}$.}
\label{fig:j54}
\end{figure}

\noindent\textit{Step 1. Hub and starting label (Figure~\ref{fig:j54-step1}).} 
Here $p=3$, $q=24$ and $b=23$. Since $m=4$ is even,
$c=pm+1=13$. The spoke set is
$S=\{23,23-6,3,3+6\}=\{3,9,17,23\}$.
Removing $S\cup\{13\}$ from $\{1,\ldots,24\}$ gives
\[
\begin{aligned}
 D=\bigl(&24,22,21,20,19,18,16,15,14,12,\\
        &11,10,8,7,6,5,4,2,1\bigr).
\end{aligned}
\]
Set $x_0=24$ and use $x_j=x_{j-1}+(-1)^j d_j$.
The sector order is $(x_1,\ldots,x_{19},x_0)$, with radial positions
$x_1,x_6,x_{11},x_{16}$. Assign $23$ to the hub and $x_0=24$
to the last outer cycle position.
\par\medskip

\noindent\textit{Step 2. First sector (Figure~\ref{fig:j54-step2}).} 
The first five differences give
\[
\begin{aligned}
 x_1&=24-24=0,&x_2&=0+22=22,\\
 x_3&=22-21=1,&x_4&=1+20=21,\\
 x_5&=21-19=2.
\end{aligned}
\]
The sector is $(0,22,1,21,2)$; its radial label $0$ gives spoke
difference $|23-0|=23$.
\par\medskip

\noindent\textit{Step 3. Second sector (Figure~\ref{fig:j54-step3}).} 
Continue with $18,16,15,14,12$:
\[
\begin{aligned}
 x_6&=2+18=20,&x_7&=20-16=4,\\
 x_8&=4+15=19,&x_9&=19-14=5,\\
 x_{10}&=5+12=17.
\end{aligned}
\]
The sector is $(20,4,19,5,17)$; its radial label $20$ gives spoke
difference $3$.
\par\medskip

\noindent\textit{Step 4. Third sector (Figure~\ref{fig:j54-step4}).} 
The next five differences give
\[
\begin{aligned}
 x_{11}&=17-11=6,&x_{12}&=6+10=16,\\
 x_{13}&=16-8=8,&x_{14}&=8+7=15,\\
 x_{15}&=15-6=9.
\end{aligned}
\]
The sector is $(6,16,8,15,9)$; its radial label $6$ gives spoke
difference $17$.
\par\medskip

\noindent\textit{Step 5. Last sector and outer cycle closure (Figure~\ref{fig:j54-step5}).} 
The remaining differences give
\[
\begin{aligned}
 x_{16}&=9+5=14,&x_{17}&=14-4=10,\\
 x_{18}&=10+2=12,&x_{19}&=12-1=11.
\end{aligned}
\]
The radial label $14$ gives spoke difference $9$.
The reserved edge joins $x_{19}=11$ to $x_0=24$, with difference
$13=c$. The completed sector order is
\[
\begin{aligned}
 &(0,22,1,21,2)\mid(20,4,19,5,17)\\
 &\mid(6,16,8,15,9)\mid(14,10,12,11,24).
\end{aligned}
\]
After rotation, the last-to-first edge joins $24$ and $0$ and has
difference $24$.
The differences in $D$, the reserved difference $13$ and the spokes
$\{3,9,17,23\}$ partition $\{1,\ldots,24\}$.
\par\medskip

\clearpage

\subsection{\texorpdfstring{$J_{6,4}$}{J(6,4)}}

\begin{figure}[p]
\centering
\makebox[\linewidth][c]{%
\begin{subfigure}[t]{0.332\linewidth}
\centering
\includegraphics[width=\linewidth]{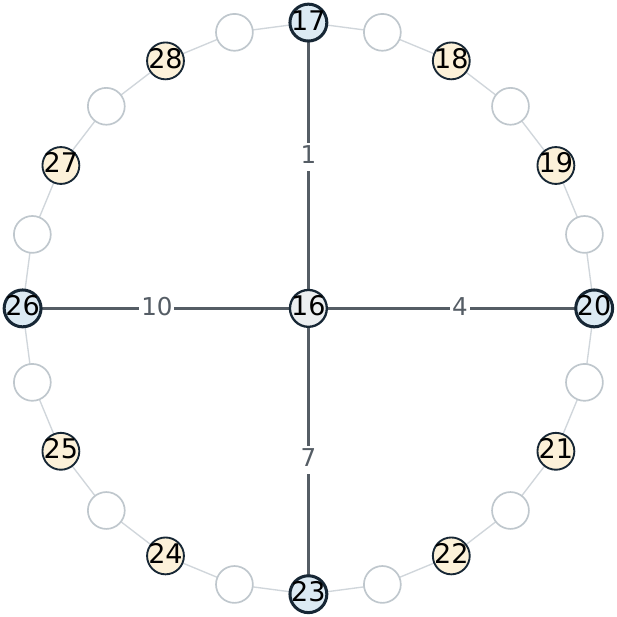}
\caption{Hub and high labels.}
\label{fig:j64-step1}
\end{subfigure}%
\hspace{\constructionpanelgap}%
\begin{subfigure}[t]{0.332\linewidth}
\centering
\includegraphics[width=\linewidth]{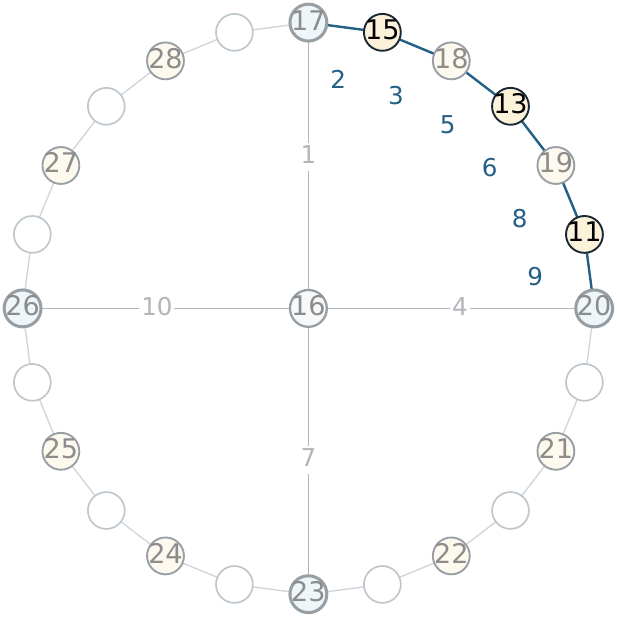}
\caption{First sector.}
\label{fig:j64-step2}
\end{subfigure}%
}\par\vspace{\constructionrowgap}
\makebox[\linewidth][c]{%
\begin{subfigure}[t]{0.332\linewidth}
\centering
\includegraphics[width=\linewidth]{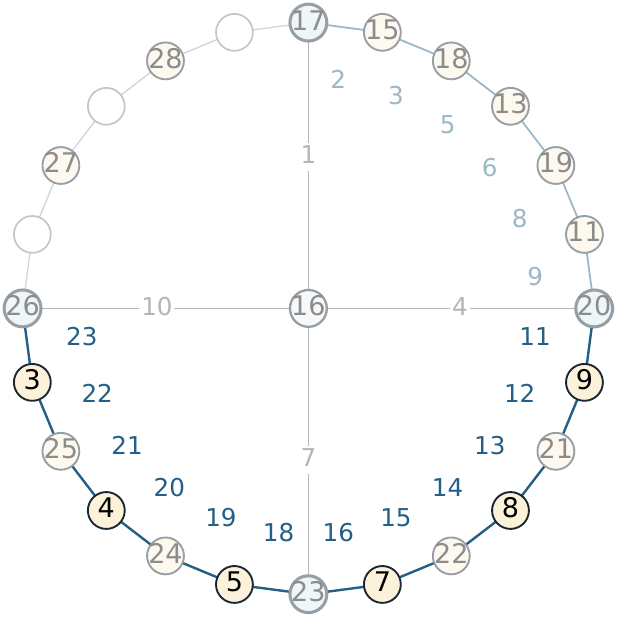}
\caption{Second and third sectors.}
\label{fig:j64-step3}
\end{subfigure}%
\hspace{\constructionpanelgap}%
\begin{subfigure}[t]{0.332\linewidth}
\centering
\includegraphics[width=\linewidth]{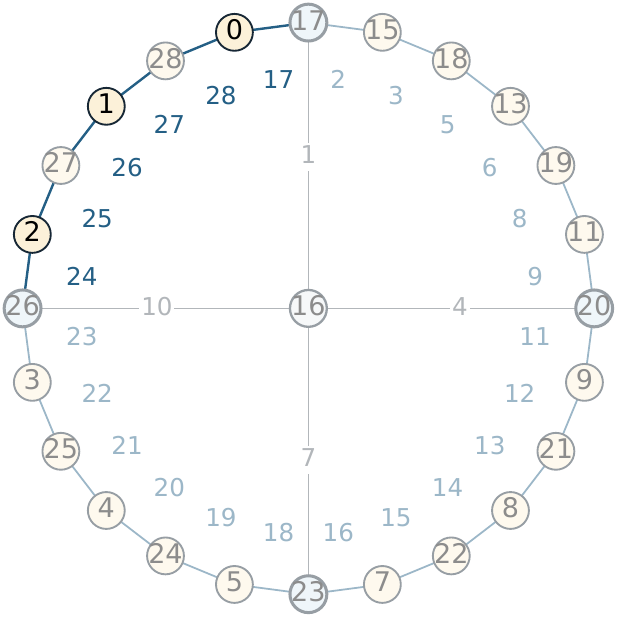}
\caption{Last sector and closure.}
\label{fig:j64-step4}
\end{subfigure}%
}\par\smallskip
\caption{Construction of $J_{6,4}$.}
\label{fig:j64}
\end{figure}

\noindent\textit{Step 1. Hub and high labels (Figure~\ref{fig:j64-step1}).} 
Here $k=3$, $P=12$, $a=16$ and $q=28$. The even-$m$ branch of
\eqref{eq:case-two} gives
\[
 b=16,\qquad \eta=0,\qquad c=17,\qquad H_j=17+j\quad(0\leq j<12).
\]
Place $H=(17,18,\ldots,28)$ at alternate outer cycle positions.
The radial labels $H_0=17$, $H_3=20$, $H_6=23$ and $H_9=26$
give $S=\{1,4,7,10\}$. Removing $S\cup\{17,28\}$ leaves
\[
\begin{aligned}
 D=\bigl(&2,3,5,6,8,9,11,12,13,14,15,16,\\
        &18,19,20,21,22,23,24,25,26,27\bigr).
\end{aligned}
\]
Every residual pair is consecutive, matching the unit increments of $H$.
\par\medskip

\noindent\textit{Step 2. First sector (Figure~\ref{fig:j64-step2}).} 
Apply $L_j=H_j-d_{2j+1}$ to the first three pairs:
\[
\begin{aligned}
 (2,3):\quad L_0&=17-2=15,\\
 (5,6):\quad L_1&=18-5=13,\\
 (8,9):\quad L_2&=19-8=11.
\end{aligned}
\]
The sector $(17,15,18,13,19,11)$ and its edge to the next radial
label $20$ realize differences $2,3,5,6,8,9$.
\par\medskip

\noindent\textit{Step 3. Second and third sectors (Figure~\ref{fig:j64-step3}).} 
The next three pairs give
\[
 L_3=20-11=9,\qquad L_4=21-13=8,\qquad L_5=22-15=7.
\]
They fill the second sector $(20,9,21,8,22,7)$.
The following pairs $(18,19)$, $(20,21)$ and $(22,23)$ give
\[
 L_6=23-18=5,\qquad L_7=24-20=4,\qquad L_8=25-22=3,
\]
filling the third sector $(23,5,24,4,25,3)$.
Low labels are at most $15$ and high labels at least $17$, leaving
$16$ for the hub.
\par\medskip

\noindent\textit{Step 4. Last sector and outer cycle closure (Figure~\ref{fig:j64-step4}).} 
The last pairs $(24,25)$ and $(26,27)$ give
\[
 L_9=26-24=2,\qquad L_{10}=27-26=1.
\]
Set $L_{11}=0$. The remaining edges have differences
$|28-0|=28$ and $|0-17|=17$. The completed outer cycle is
\[
\begin{aligned}
 &(17,15,18,13,19,11)\mid(20,9,21,8,22,7)\\
 &\mid(23,5,24,4,25,3)\mid(26,2,27,1,28,0).
\end{aligned}
\]
Its outer cycle differences are the values in $D$ together with $17$ and $28$; the spokes supply
$\{1,4,7,10\}$. Together they give $\{1,\ldots,28\}$ exactly once.
\par\medskip

\clearpage

\subsection{\texorpdfstring{$J_{7,4}$}{J(7,4)}}

\begin{figure}[p]
\centering
\makebox[\linewidth][c]{%
\begin{subfigure}[t]{0.332\linewidth}
\centering
\includegraphics[width=\linewidth]{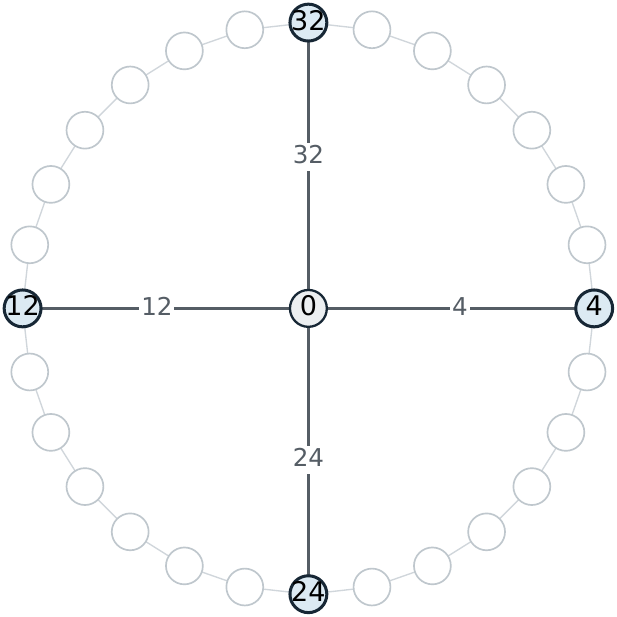}
\caption{Hub and radial labels.}
\label{fig:j74-step1}
\end{subfigure}%
\hspace{\constructionpanelgap}%
\begin{subfigure}[t]{0.332\linewidth}
\centering
\includegraphics[width=\linewidth]{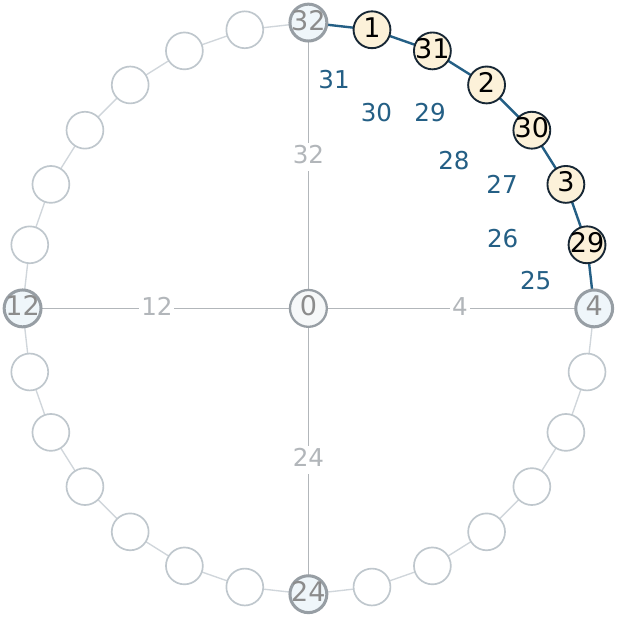}
\caption{First sector.}
\label{fig:j74-step2}
\end{subfigure}%
\hspace{\constructionpanelgap}%
\begin{subfigure}[t]{0.332\linewidth}
\centering
\includegraphics[width=\linewidth]{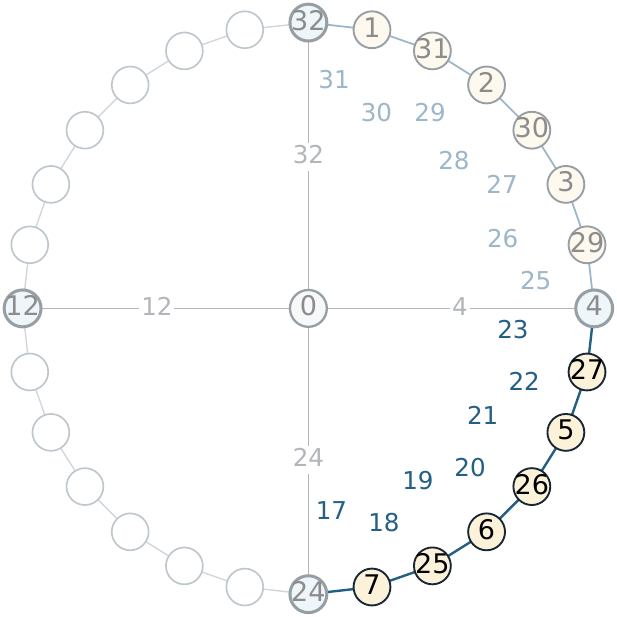}
\caption{Second sector.}
\label{fig:j74-step3}
\end{subfigure}%
}\par\vspace{\constructionrowgap}
\makebox[\linewidth][c]{%
\begin{subfigure}[t]{0.332\linewidth}
\centering
\includegraphics[width=\linewidth]{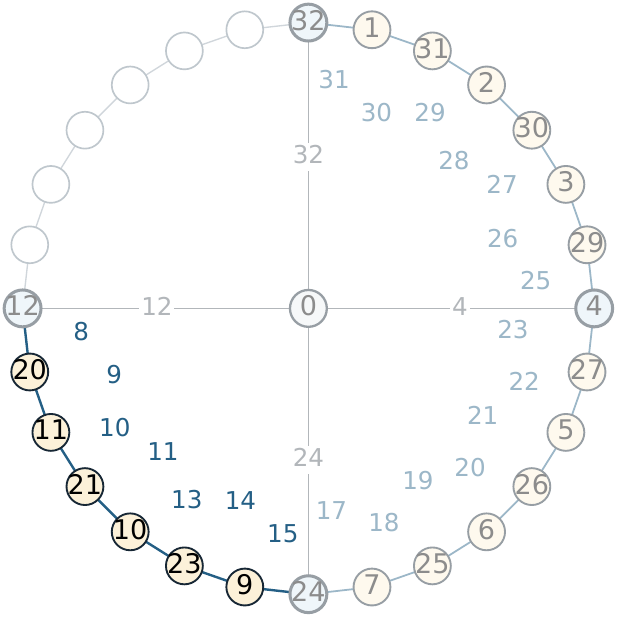}
\caption{Third sector.}
\label{fig:j74-step4}
\end{subfigure}%
\hspace{\constructionpanelgap}%
\begin{subfigure}[t]{0.332\linewidth}
\centering
\includegraphics[width=\linewidth]{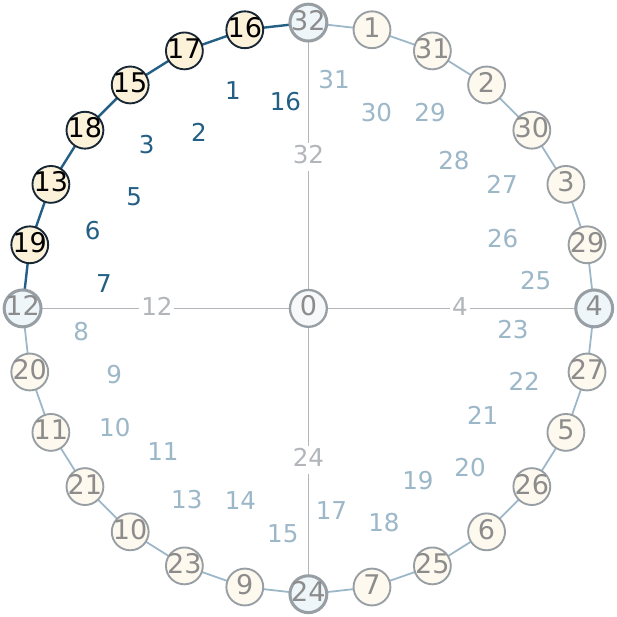}
\caption{Last sector and closure.}
\label{fig:j74-step5}
\end{subfigure}%
}\par\smallskip
\caption{Construction of $J_{7,4}$.}
\label{fig:j74}
\end{figure}

\noindent\textit{Step 1. Hub and radial labels (Figure~\ref{fig:j74-step1}).} 
Here $p=4$, $r=2$, $M=2$ and $q=32$. The sector parameters are
\[
 A=(0,4,8,12),\qquad B=(32,28,24,20),\qquad \delta=(0,0,1,1).
\]
Assign $0$ to the hub. The radial labels are $B_0=32$, $A_1=4$,
$B_2=24$ and $A_3=12$, so $S=\{4,12,24,32\}$.
Reserve $c=pm=16$. The other outer cycle differences, in decreasing order, are
\[
\begin{aligned}
 D=\bigl(&31,30,29,28,27,26,25,23,22,21,20,19,18,17,\\
        &15,14,13,11,10,9,8,7,6,5,3,2,1\bigr).
\end{aligned}
\]
\par\medskip

\noindent\textit{Step 2. First sector (Figure~\ref{fig:j74-step2}).} 
For $i=0$ and $\delta_0=0$, \eqref{eq:case-three} gives
\[
\begin{aligned}
 \lab(v_{0,2s})&=32-s &&(0\leq s\leq3),\\
 \lab(v_{0,2s+1})&=s+1 &&(0\leq s\leq2).
\end{aligned}
\]
The sector is $(32,1,31,2,30,3,29)$.
Its internal differences are $31,30,29,28,27,26$; the edge to the next
radial label $4$ has difference $25$.
\par\medskip

\noindent\textit{Step 3. Second sector (Figure~\ref{fig:j74-step3}).} 
For $i=1$, $\delta_1=0$, so
\[
\begin{aligned}
 \lab(v_{1,2s})&=4+s &&(0\leq s\leq3),\\
 \lab(v_{1,2s+1})&=27-s &&(0\leq s\leq2).
\end{aligned}
\]
The sector is $(4,27,5,26,6,25,7)$.
Its internal differences are $23,22,21,20,19,18$; the edge to the next
radial label $24$ has difference $17$.
\par\medskip

\noindent\textit{Step 4. Third sector (Figure~\ref{fig:j74-step4}).} 
For $i=2$, $\delta_2=1$ activates the correction:
\[
\begin{aligned}
 \lab(v_{2,2s})&=24-s-\ind{s\geq2} &&(0\leq s\leq3),\\
 \lab(v_{2,2s+1})&=9+s &&(0\leq s\leq2).
\end{aligned}
\]
The even positions receive $24,23,21,20$ and the odd positions
receive $9,10,11$. In particular, $\lab(v_{2,4})=24-2-1=21$.
The sector $(24,9,23,10,21,11,20)$ has internal differences
$15,14,13,11,10,9$. Its edge to radial label $12$ has difference $8$.
\par\medskip

\noindent\textit{Step 5. Last sector and outer cycle closure (Figure~\ref{fig:j74-step5}).} 
For $i=3$, $\delta_3=1$, giving
\[
\begin{aligned}
 \lab(v_{3,2s})&=12+s+\ind{s\geq2} &&(0\leq s\leq3),\\
 \lab(v_{3,2s+1})&=19-s &&(0\leq s\leq2).
\end{aligned}
\]
The sector is $(12,19,13,18,15,17,16)$, with internal differences
$7,6,5,3,2,1$. The closing difference is $|16-32|=16$.
The completed outer cycle is
\[
\begin{aligned}
 &(32,1,31,2,30,3,29)\mid(4,27,5,26,6,25,7)\\
 &\mid(24,9,23,10,21,11,20)\mid(12,19,13,18,15,17,16).
\end{aligned}
\]
The differences in $D$, the closing difference $16$ and the spoke
differences $\{4,12,24,32\}$ partition $\{1,\ldots,32\}$.
\par\medskip

\clearpage

\section*{Acknowledgement}
The research that led to the present paper was partially supported by a grant
of the group of GNSAGA of INdAM, Italy.

\Needspace{12\baselineskip}


\begin{thebibliography}{99}
\setlength{\itemsep}{2pt}
\raggedright

\bibitem{BondyMurty2008}
J.~A. Bondy and U.~S.~R. Murty. \emph{Graph Theory}.
Graduate Texts in Mathematics, volume 244. Springer, 2008.
\doi{10.1007/978-1-84628-970-5}.

\bibitem{EdwardsHoward2006}
M.~Edwards and L.~Howard. A survey of graceful trees.
\emph{Atlantic Electronic Journal of Mathematics}, 1(1):5--30, 2006.
\newline\url{http://euclid.trentu.ca/aejm/V1N1/Edwards.v1n1.pdf}.

\bibitem{ElZanatiVandenEynden2009}
S.~I. El-Zanati and C.~Vanden Eynden.
On Rosa-type labelings and cyclic graph decompositions.
\emph{Mathematica Slovaca}, 59(1):1--18, 2009.
\doi{10.2478/s12175-008-0108-x}.

\bibitem{Gallian2025}
J.~A. Gallian. A dynamic survey of graph labeling.
\emph{The Electronic Journal of Combinatorics} 5, 1998.
DS6, twenty-eighth edition, October 30, 2025. \doi{10.37236/27}.

\bibitem{Gardner1972}
M.~Gardner. Mathematical games: The graceful graphs of Solomon Golomb,
or how to number a graph parsimoniously.
\emph{Scientific American}, 226(3):108--112, March 1972.
\doi{10.1038/scientificamerican0372-108}.

\bibitem{Golomb2008}
S.~W. Golomb. Golomb's puzzle column: Graceful graphs solutions.
\emph{IEEE Information Theory Society Newsletter}, 58(3):17--18,
September 2008.
\newline\url{https://www.itsoc.org/sites/default/files/2021-03/itNL0908.pdf}.

\bibitem{GrahamSloane1980}
R.~L. Graham and N.~J.~A. Sloane. On additive bases and harmonious graphs.
\emph{SIAM Journal on Algebraic and Discrete Methods}, 1(4):382--404, 1980.
\doi{10.1137/0601045}.

\bibitem{HaviarIvan2026}
M.~Haviar and J.~Ivan. Characterizations of kites as graceful graphs II.
\emph{CUBO, A Mathematical Journal}, 28(3):523--544, 2026.
\doi{10.56754/0719-0646.2803.523}.

\bibitem{HaviarKotulova2024}
M.~Haviar and K.~Kotu\v{l}ov\'a.
Characterizations of kites as graceful graphs.
\emph{CUBO, A Mathematical Journal}, 26(3):367--386, 2024.
\doi{10.56754/0719-0646.2603.367}.

\bibitem{HaviarKurtulik2021}
M.~Haviar and S.~Kurtul\'ik.
A new class of graceful graphs: $k$-enriched fan graphs and their characterisations.
\emph{CUBO, A Mathematical Journal}, 23(2):313--331, 2021.
\doi{10.4067/S0719-06462021000200313}.

\bibitem{ImbesiLaBarbiera2024}
M.~Imbesi and M.~La Barbiera.
On the Betti numbers and gracefulness of some planar graphs.
\arxiv{1210.8198v2}, March 29, 2024.

\bibitem{LopezMuntaner2017}
S.~C. L\'opez and F.~A. Muntaner-Batle.
\emph{Graceful, Harmonious and Magic Type Labelings: Relations and Techniques}.
SpringerBriefs in Mathematics. Springer, 2017.
\doi{10.1007/978-3-319-52657-7}.

\bibitem{Rofa2021}
R.~I. Rofa. A graceful algebraic function labelling of rooted symmetric trees.
\arxiv{2109.09511}, September 15, 2021.

\bibitem{Rosa1967}
A.~Rosa. On certain valuations of the vertices of a graph.
In \emph{Theory of Graphs: International Symposium, Rome, July 1966},
pages 349--355. Gordon and Breach and Dunod, 1967.
\href{https://www.cs.columbia.edu/~tim/teaching/cs4203/Rosa-GracefulLabelings.pdf}{Online scan}.

\end{thebibliography}
\end{document}